\documentclass[12pt]{article}

\usepackage{fullpage}
\usepackage{multicol,multirow}
\usepackage{tabularx}
\usepackage{ulem}
\usepackage[T2A]{fontenc}
\usepackage[utf8]{inputenc}
\usepackage[russian]{babel}
\usepackage{amssymb}
\usepackage{amsfonts}
\usepackage{latexsym}
\usepackage{amsmath}
\usepackage{color}

\begin{document}

\newcommand{\Expect}{\mathsf{E}}
\newtheorem{theorem}{Theorem}[section]
\newtheorem{definition}{Definition}[section]
\newtheorem{lemma}{Lemma}[section]
\newtheorem{remark}{Remark}[section]

\begin{center}
\LARGE{\textbf{Branching Random Walks with Immigration.\\Lyapunov Stability}}
\end{center}

\begin{center}
Yu. Makarova$^1$, D. Han$^2$, S. Molchanov$^{3,4}$, E. Yarovaya$^1$
\end{center}

\begin{center}
$^1$ Lomonosov Moscow State University,\\
Leninskie Gory 1, Moscow, 119234, Russia\\
\textit{ykmakarova@gmail.com, yarovaya@mech.math.msu.su}\\
$^2$ University of Louisville,\\
Louisville, KY 440292, USA\\
\textit{mathbaobao@gmail.com}\\
$^3$ University of North Carolina at Charlotte,\\
Charlotte, NC 28223, USA\\
$^4$ National Research University Higher School of Economics,\\
Mysnitskaya str., Moscow, 101000, Russia\\
\textit{smolchan@uncc.edu}
\end{center}   

\section{Introduction}
\label{introduction}

\quad The initial version of this article was published in proceedings of the international scientific conference ACMPT 2017 which was dedicated to the 90th birth anniversary of outstanding mathematician and expert of applied mathematics Professor A.D.Solov'ev. The extended version was prepared at that time. But it was not published due to some technical difficulties. The present article is the final version of the initial one. 

Nowadays the branching random walks (BRWs) are an appropriate tool to describe and explore the evolution processes with birth, death and migration~\cite{zeldaoth}. In practice such models may be used in biology~\cite{emy} and demography~\cite{MW2017}.

Models presented in the paper give a reasonably good description for the demographic situation associated with immigration in different Europian countries. Population dynamics research the statistical equilibrium of the process, so-called \textit{steady state}. One of the examples of such stochastic processes is a continuous-time critical Galton-Watson branching process where the rates of birth and annihilation are equal. This process with a random walk of particles and their  generation at any point on the lattice was considered in~\cite{HMW2017}. Such BRW has steady state under some conditions. In~\cite{HMW2017} authors assume  that the underlying random walk associated with the process is transient. Unfortunately, this population model is not stable with respect to (even small) random perturbation and the steady state destroys, see~\cite{KKP:IDAQP08}. It is not difficult to understand the fact: in case when death rate $\mu$ and birth rate $\beta$ in model of binary splitting are equal and the perturbation has the form $\beta' = \beta + \varepsilon$, $\mu'=\mu$ (where $\beta'$ and $\mu'$ new rates of birth and death), then the process becomes supercritical for any small $\varepsilon > 0$. If the perturbation has the form $\beta' = \beta-\varepsilon$, then the process becomes subcritical and the population degenerates.

Typically, immigration was introduced in models of branching processes. One of the first models was introduced by B.~Sevastyanov in~\cite{sevast}. He considered the model when each particle can produce an arbitrary number of offsprings and any number of particles could appear from the outside in the system. We are going to consider immigration in more complex BRW model, in which also containes a random walk in space.

The arising problem of absence the steady state can be solved by adding immigration. The presence of immigration can stabialize the process and stop extinction when the birth rate is less than the death rate. Apparently, this approach was suggested by Han, Molchanov and Whitmeyer in~\cite{HMW2017}, but only for the case of binary splitting, thus is each particle can produce only one offspring. For such BRW it is possible to have slightly different interpretation: each particle produce two offsprings which start their evolution processes independently and the parental particle die.

The structure of the paper is following. In section~\ref{description_of_the_model} we describe the model of the BRW with immigration the main feature of which is arbitrary splitting of particles. Moreover, we present some conditional expectations. In section~\ref{first_moment} we study the first moment of the particle field, derive the differentional equation and obtain the asymptotic behaviour of the solution. In section~\ref{second_moment} we consider the second moment of the particle field. Similarly to the previous section, we get the asymptotic behaviour of the second moment, which can be obtained from the derived differentional equation. Unlike the first moment, derivation of the equation of the second moment needs the usage of more complicated tools. Besides, we should research two cases: they depend on the equality of points on the lattice which we consider. In section~\ref{n_moment} there are the differentional equations of the higher moments with their derivations. In section~\ref{gen_function} we consider the generating function of the process. Typically, this is a useful tool in studying the process, because it gives an opportunity to get all the moments of the random variable. However in our case is not a good method for research due to the fact that the usage of generating function does not allow to obtain the equations in all necessary cases. In section~\ref{lyapunov_stab} we assume that the intensities of the birth, death and immigration are functions which depend on the points on the lattice. Under these assumptions we consider the Lyapunov stability. This is more realistic in applications~\cite{MW2017}. For example, people prefer to live in areas where there are more resources required for living, rather than in areas with uncomfortable climate. 

\section{BRW model on the Multidimensional Lattice}
\label{description_of_the_model}

\quad In this section we consider countinious-time symmetric BRW on the lattice $\mathbb Z^d$. The subject of the study is the particle field $n(t,x)$, where $t \geqslant 0$, $x \in \mathbb Z^d$. We assume that at the initial moment random variables $n(0,x)$ are independent and identically distributed with finite moments. For example, it can be said the $n(0,x) \equiv 1$ for any $x\in\mathbb Z^d$. In fact, we will show in the future that for $t\to\infty$ the contribution of $n(0,x)$ into $n(t,x)$ is exponentionally neglectable. The evolution of the particle field includes several options.

Each particle at the moment$t > 0$ in the point $x \in \mathbb Z^d$ stays at this point random time $\tau$ up to the first transformation. The process we consider is the Markov process, so the random variable $\tau$ has the exponentional distribution. Therefore, at the moment $t + \tau + 0$ there can be the fololowing transformations:
\begin{enumerate}
\item Firstly, it can be the jump from the point $x$ to the point $x+z$ with the probability $a(z)$. We assume that $a(z) = a(-z)$, $\sum_{z \ne 0} a(z) = 1$ and $a(0) = -1$.  The intensity of jump, also called the diffusion coefficient, is denoted by $\varkappa > 0$. Thus, the probability to jump from the point $x$ to the point $x+z$ during the small time $dt$ is $\varkappa a(z) dt$. The generator of underlying random form has the form:
\begin{align*}
(\mathcal L \psi)(x) = \varkappa \sum_{z \ne 0} [\psi(x+z) - \psi(x)]a(z)\enspace.
\end{align*}

Moreover, we assume that our random walk is irreducible, thus is $Span\{ z: a(z) > 0 \} = \mathbb Z^d$.

The operator $\mathcal L$ generates the Markov semigroup
\begin{align*}
P_t = \exp\{t\mathcal L\}\enspace.
\end{align*}
with kernel (transition probability) $p(t,x,y)$. Here
\begin{align*}
\displaystyle\frac{\partial p(t,x,y)}{\partial t} = \mathcal L_x p(t,x,y) = \mathcal L_y p(t,x,y)\enspace.
\end{align*}

Each particle carries out random walk with the generator $\mathcal L$ up to the moment of the first transformation, which consists of either annihilation or splitting.

\item Secondly, particle can die with the probability $\mu dt$, where $\mu$ is the mortality rate.

\item Thirdly, each particle (independent on others) can produce $n$ offsprings (or we can say that it produces $n-1$ particles and still stays at the same point on the lattice). Let $b_n$, $n \ne 1$ is the intencity of the transformation for the single parental particle into $n$ particles. Besides,
\begin{align*}
\mu + \sum_{n\geqslant 2} b_n = -b_1 > 0\enspace.
\end{align*}

So we obtain the infinitesimal generating function:
\begin{align*}
F(z) = \mu + \sum_{n=1}^{\infty} b_n z^n = \mu - (\mu + \sum_{n\geqslant 2} b_n)z + \sum_{n\geqslant 2} b_n z^n\enspace.
\end{align*}

We also assume that $F(z)$ is the analytical function in the circle $|z| < 1 + \delta$, $\delta > 0$, thus is the intensities $b_n$ as the functions of $n$ are exponentially decreasing. Lastly, new particles start their evolution independently on others at the points where they appear.

\item Finally, the new property of the process is the presence of immigration. It means that a new particle can appear at the time interval $(t,t+dt)$ at any point $x\in\mathbb Z^d$ on the lattice with the probability $kdt$, where $k$ is the immigration rate.

\end{enumerate}

Futher we assume that the intensities $b_n$, $\mu$ and $k$ are constant. 

The study of the BRW is usually based on backward Kolmogorov equations:
\begin{align*}
\dot P = \mathcal A P\enspace,
\end{align*}
where $\mathcal A$ is the transition probability matrix.

However in our model, because of the presence of immigration and random distribution of the initial particle field, we have to use forward Kolmogorov equations:
\begin{align*}
\dot P = P \mathcal A\enspace.
\end{align*}

The derivation of the forward Kolmogorov equation in our model is based on the following representation:
\begin{align*}
n(t+dt,x) = n(t,x) + \xi(dt,x)\enspace,
\end{align*}
where $\xi(dt,x)$ is the discrete random variable with the distribution
\begin{equation*}
\xi (dt,x) = \begin{cases}
n-1, &\text{with probability}~b_{n} n(t,x)dt,~ n \geqslant 3,\\
1, &\text{with probability}~b_{2} n(t,x)dt + kdt\\
&+ \varkappa \sum_{z \ne 0} a(-z)n(t,x+z)dt,\\
-1, &\text{with probability}~\mu n(t,x)dt + \varkappa n(t,x)dt,\\
0, &\text{with probability}~1 - \sum_{n \geqslant 3} b_{n} n(t,x)dt\\
&- (b_2 + \mu + \varkappa)n(t,x)dt - kdt\\
& -\sum_{z \ne 0} a(-z)n(t,x+z)dt.
\end{cases}
\end{equation*}

One of the targets of our research is to find the asymptotic behaviour of the moments of the random variable $n(t,x)$ (definition of moments will be in the next sections). In this case we need the technique of conditional expectations\medspace\cite[гл.~2]{shiryaev}.

Notice that random variable $\xi(dt,x)$ and the $\sigma$-algebra $\mathcal F_{\leqslant t}$, where $\mathcal F_{\leqslant t}$ is the $\sigma$-algebra of events before and including $t$, are independent.

In the end of the section we consider some useful relations.

1.Firstly, we have the conditional expectation of the variable $\xi(dt,x)$.
\begin{multline}\label{first_ce}
\Expect[ \xi(dt,x)| \mathcal F_{\leqslant t}] = \sum_{n=2}^{\infty} (n-1) b_n n(t,x)dt + kdt + \sum_{z \ne 0} \varkappa a(-z)n(t,x+z)dt \\- (\mu + \varkappa)n(t,x)dt;
\end{multline}

2. The expectation of $\xi^2 (dt,x)$ is
\begin{multline}\label{second_ce}
\Expect[ \xi^2(dt,x)| \mathcal F_{\leqslant t}] = \sum_{n=2}^{\infty} (n-1)^2 b_n n(t,x)dt + kdt + \sum_{z \ne 0} \varkappa a(-z)n(t,x+z)dt \\+ (\mu + \varkappa)n(t,x)dt;
\end{multline}

3. The correlation of $\xi(dt,x)$ and $\xi(dt,y)$. Here we assume that $dt^2 = 0$.
\begin{multline}\label{third_ce}
\Expect[ \xi(dt,x) \xi(dt,y)| \mathcal F_{\leqslant t}] = -\varkappa (a(y-x)n(t,x)dt + a(x-y)n(t,y)dt),\medspace x \ne y;
\end{multline}

4. For three different points $x$, $y$, $z$ $\in\mathbb Z^d$ (using $dt^2=0$) we get
\begin{multline}\label{fourth_ce}
\Expect[ \xi(dt,x) \xi(dt,y) \xi(dt,z)| \mathcal F_{\leqslant t}] = 0,\medspace x \ne y,\medspace x \ne z,\medspace x \ne z;
\end{multline}

5. For the higher moments of the particle field we have the expectation of $\xi^p (dt,x)$, $p\geqslant 2$.
\begin{multline}\label{seventh_ce}
\Expect[ \xi(dt,x)^{p}| \mathcal F_{\leqslant t}] = \sum_{n=2}^{\infty} (n-1)^p b_ndt + kdt + \varkappa \sum_{z \ne 0} a(z)n(t,x+z)dt \\+ (-1)^{p}[\mu + \varkappa]n(t,x)dt,\medspace p>0;
\end{multline}

6. The next relation is also useful for getting the equations for the higher moments.
\begin{multline}\label{eighth_ce}
\Expect[ \xi(dt,x)^{p} \xi(dt,y)^{q}| \mathcal F_{\leqslant t}] = \varkappa [(-1)^{p}a(y-x)n(t,x)dt \\+ (-1)^{q}a(x-y)n(t,y)dt],\medspace x \ne y,\medspace p,q>0;
\end{multline}

7. Finally, we present the relation which will be used in method of generationg fuction.
\begin{multline}\label{nineth_ce}
\Expect[ e^{-z\xi(dt,x)}| \mathcal F_{\leqslant t}] = \sum_{n\geqslant2} e^{-z(n-1)}b_n n(t,x)dt \\+ e^{-z}(kdt + \varkappa \sum_{y\ne0} a(y)n(t,x+y)dt) + e^z (\mu + \varkappa)n(t,x)dt \\+ (1 - \sum_{n\geqslant2} b_n n(t,x)dt - kdt - \varkappa \sum_{y\ne0} a(y) n(t,x+y)dt \\- \mu n(t,x)dt - \varkappa n(t,x)dt),\medspace z \geqslant 0.
\end{multline}

\section{The First Moment}
\label{first_moment}

\quad In the future study there will be used some properties of the conditional expectations, see~\cite{shiryaev} for details:
\begin{equation}\label{indep}
\text{If } \xi - \mathcal G-\text{measurable, then } \Expect(\xi\eta|\mathcal G) = \xi \Expect(\eta| \mathcal G)\enspace,
\end{equation}
\begin{equation}\label{fullprob}
\Expect(\Expect(\xi|\mathcal G)) = \Expect\xi\enspace.
\end{equation}

Define the first moment for the particle field $n(t,x)$ as following:
\begin{align*}
m_1(t,x) = \Expect n(t,x)\enspace.
\end{align*}

Now we are going to get the differentional equation for the first moment. Consider the first moment at the time $t+dt$ to obtain the eqauation:
\begin{align*}
En(t+dt,x) &= E[E[n(t+dt,x)|\mathcal F_{\leqslant t}]] = E[E[n(t,x) + \xi (dt,x)|\mathcal F_{\leqslant t}]]\enspace.
\end{align*}

Using properties~(\ref{indep}) and~(\ref{fullprob}) and equality~(\ref{first_ce}) we have

\begin{align*}
En(t+dt,x) &= E[E[n(t+dt,x)|\mathcal F_{\leqslant t}]] = E[E[n(t,x) + \xi (dt,x)|\mathcal F_{\leqslant t}]] \\&= m_1 (t,x) + \sum_{n=2}^{\infty} (n-1) b_n m_1 (t,x)dt + kdt \\&+ \sum_{z \ne 0} \varkappa a(z)(m_1 (t,x+z)dt  - m_1 (t,x)dt - \mu m_1 (t,x)dt\enspace.
\end{align*}

Let
\begin{align*}
\mathcal L_a f(t,x) := \sum_{z \ne 0} a(z)\Bigl(f(t, x+z) - f(t,x)\Bigr)\enspace.
\end{align*}

Combining all abobe results, we get the differentional equation for the first moment. Moreover, due to the space homogenety, $\mathcal L_a m_1(t,x) = 0$.
\begin{equation}\label{diff_eq_f_mom_const_coef}
\begin{cases}
\displaystyle\frac {\partial{m_1(t,x)}}{\partial{t}} &= \Bigl(\sum_{n=2}^{\infty} (n-1) b_n - \mu\Bigr) m_1 (t,x) + k\enspace,\\[4pt]
m_1(0,x) &= En(0,x)\enspace.
\end{cases}
\end{equation}

Let us define the coefficient $\beta$ --- the birth rate as
\begin{align*}
\beta := \sum_{n\geqslant 2} (n-1) b_n\enspace.
\end{align*}

Equation~(\ref{diff_eq_f_mom_const_coef}) is similar to equation which covers the case when $\beta = \beta(x)$, $\mu = \mu(x)$, $k = k(x)$ are bounded functions on the lattice $\mathbb Z^d$. This case will be considered later.

In case of constant coefficients the equation can be solved as an ordinary differential equation. The solution has the form
\begin{align*}
m_1 (t,x) = \frac {k}{\beta - \mu} (e^{(\beta - \mu)t} - 1) + e^{(\beta- \mu)t} En(0,x)\enspace.
\end{align*}

\begin{remark}
If $\beta > \mu$ (in this case we say that the process is supercritical) then the population  grows exponentially, despite the presence of immigration.

If $\beta = \mu$ (such case is critical) and $k > 0$ then  the population grow with linear speed.

We are interested in the last case (subcritical case) when $\mu > \beta$. Here, if $k = 0$ then the population vanishes. But if $k > 0$ then
\begin{align}\label{as_first_mom}
m_1(t,x) \to \frac{k}{\mu - \beta},\medspace t \to \infty\enspace.
\end{align}
\end{remark}

Futher we are going to consider the case when $\mu > \beta$ and $k > 0$.

In case of non-constant coefficients we obtain the Lyapunov stability for the first moment. See the section~\ref{lyapunov_stab} for details.

\begin{theorem}
Let $b_n (x)$, $n \geqslant 2$, $\mu (x)$, $k(x)$, $x \in \mathbb Z^d$ are bounded and $\mu (x) - \beta (x) \geqslant \delta_1 > 0$, $k(x) \geqslant \delta_2 > 0$. Then for the bounded initial conditions there exists the limit
\begin{equation*}
m_1 (\infty,x) = \lim_{t \to \infty} m_1 (t,x)\enspace.
\end{equation*}
\end{theorem}

\section{The Second Moment}
\label{second_moment}

\quad Let us denote the second moment:
\begin{align*}
m_2(t,x,y) = \Expect [n(t,x)n(t,y)]\enspace.
\end{align*}

To find the asymptotic behaviour of the second moment we calculate the differential equations in two cases ($x = y$ and $x \ne y$) and combine them into one equation, using the properties of the BRW and relations~(\ref{as_first_mom}),~\eqref{indep} and~\eqref{fullprob} obtained for the first moment.

\subsection{Case 1. $x = y$}
\label{sec_mom_case_1}

\quad Here there is used the same technique as in section~\ref{first_moment} when the equation for the first moment was received. It means that we consider the second moment at the time $t+dt$ and use properties~(\ref{indep}) and~(\ref{fullprob}) and equalities~(\ref{first_ce}),~(\ref{second_ce}),~(\ref{third_ce}). To simplify the recording, we introduce the designation
\begin{align*}
\mathcal L_{ax} f(t,x,y) := \sum_{z \ne 0} a(z)\Bigl(f(t,x+z,y) - f(t,x,y)\Bigr)\enspace.
\end{align*}

Then
\begin{align*}
m_2 (t+dt,x,x) &= En^2 (t+dt,x) = E[E[n^2 (t+dt,x)|\mathcal F_{\leqslant t}]] \\& = E[E[(n(t,x) + \xi (dt,x))^2|\mathcal F_{\leqslant t}]] \\& = 2 ( \beta - \mu ) m_2 (t,x,x)dt + 2 \varkappa \mathcal L_{ax} m_2 (t,x,x)dt \\& + 2k m_1(t,x)dt + kdt + \sum_{n=2}^{\infty} (n-1)^2 b_n m_1(t,x)dt \\& + \varkappa \mathcal L_a m_1(t,x)dt + 2 \varkappa m_1 (t,x)dt + \mu m_1 (t,x)dt\enspace.
\end{align*}

From this and~\eqref{as_first_mom} it is easy to get the differential equation in case 1:
\begin{equation}\label{sec_mom_case_1_eq}
\begin{cases}
\displaystyle\frac {\partial{m_2 (t,x,x)}}{\partial{t}} &= 2 m_2 (t,x,x)[ \sum_{n=2}^{\infty} (n-1) b_n - \mu] \\[4pt]&+ \frac {k\bigl(2k + 2 \varkappa + 2 \mu + \sum_{n=2}^{\infty} (n-1)(n-2) b_n\bigr)}{\mu - \sum_{n=2}^{\infty} (n-1) b_n} \\[4pt]
&+ 2 \varkappa \mathcal L_{ax} m_2(t,x,x)\enspace,\\[4pt]
m_2(0,x,x) &= En^2(0,x)\enspace.
\end{cases}
\end{equation}

\subsection{Case 2. $x \ne y$}
\label{sec_mom_case_2}

\quad As in~\ref{sec_mom_case_1} explore the second moment at the time $t+dt$. Let
\begin{align*}
\mathcal L_{ay} f(t,x,y) := \sum_{z \ne 0} a(z)\Bigl(f(t,x,y+z) - f(t,x,y)\Bigr)\enspace.
\end{align*}

In case $x \ne y$ the folowing representation is true
\begin{align*}
m_2 (t+dt,x,y) &= E[E[n(t+dt,x)n(t+dt,y)| \mathcal F_{\leqslant t}]] \\& = E[E[(n(t,x) + \xi (dt,x))(n(t,y) + \xi (dt,y))| \mathcal F_{\leqslant t}]] \\& = m_2(t,x,y) + m_2(t,x,y)(2 \beta - 2 \mu)dt + \varkappa \mathcal L_{ax} m_2(t,x,y) \\&+ \varkappa \mathcal L_{ay} m_2 (t,x,y) + k\Bigl(m_1 (t,y) + m_1 (t,x)\Bigr)dt \\&- \varkappa\Bigl(a(y-x) m_1 (t,x) + a(x-y) m_1 (t,y)\Bigr)dt\enspace.
\end{align*}

So the differentional equation for the second case has the form:
\begin{equation}\label{sec_mom_case_2_eq}
\begin{cases}
\displaystyle\frac {\partial{m_2(t,x,y)}}{\partial{t}} &= m_2(t,x,y)(2 \beta - 2 \mu) + \varkappa \mathcal L_{ax} m_2(t,x,y) \\[4pt]&
+ \varkappa \mathcal L_{ay} m_2 (t,x,y) + k \bigl(m_1 (t,x) + m_1 (t,y)\bigr) \\[4pt]&- \varkappa\bigl(a(y-x)m_1 (t,x) + a(y-x)m_1 (t,y)\bigr)\enspace,\\[4pt]
m_2 (0,x,y) &= (En(0,x))^2\enspace.
\end{cases}
\end{equation}

\subsection{Differentional Equation for the Second Moment}
\label{eq_sec_mom}

\quad To obtain the differentional equation for the second moment we should combine equations in~\ref{sec_mom_case_1} and~\ref{sec_mom_case_2}. Note that for fixed $t$ the number of particles $n(t,x)$ homogeneous in space, therefore, it is possible to write
\begin{align*}
m_2(t,x,y) = m_2(t,x-y) = m_2(t,u)\enspace.
\end{align*}

Thus, the equation which combines~(\ref{sec_mom_case_1_eq}) and~(\ref{sec_mom_case_2_eq}) is
\begin{equation*}
\begin{cases}
\displaystyle\frac {\partial{m_2(t,u)}}{\partial{t}} &= 2 m_2(t,u)\bigl(\beta - \mu\bigr) + 2\varkappa \mathcal L_{au} m_2 (t,u) + 2\varkappa a(u) \Phi(m_1) \\[4pt]&+ \delta_0 (u) \Psi(m_1)\enspace,\\[4pt]
m_2 (0,u) &= (En(0,u))^2 (1 - \delta_0 (u)) + \delta_0 (u) En^2(0,u)\enspace.
\end{cases}
\end{equation*}

Here $x-y=u$, functions $\Phi(x)$ and $\Psi(x)$ are know functions which depend linearly on $x$.

The result obtained in~\eqref{as_first_mom} (the relation fo the first moment of the particle field) allows to write the final differentional equation for the second moment:
\begin{equation}\label{sec_mom_diff_eq}
\begin{cases}
\displaystyle\frac {\partial{m_2(t,u)}}{\partial{t}} &= 2 m_2(t,u)\bigl(\beta - \mu\bigr) + 2\varkappa \mathcal L_{au} m_2 (t,u) + \frac {2k^2}{\mu - \beta} - \frac{2\varkappa k a(u)}{\mu - \beta}\\[4pt]&\quad + \delta_0 (u) \frac{k(2\mu + \sum_{n \geqslant 2} (n-1)(n-2) b_n)}{\mu - \beta}\enspace, \\[4pt]
m_2 (0,u) &= (En(0,u))^2 (1 - \delta_0 (u)) + \delta_0 (u) En^2(0,u)\enspace.
\end{cases}
\end{equation}

\subsection{Asymptotic behaviour of the Second Moment}

\quad The next goal is to solve the equation~(\ref{sec_mom_diff_eq}) and find the asymptotic behaviour of the second moment when $t \to \infty$. One way to resolve this problem is to divide~(\ref{sec_mom_diff_eq}) into three equations and solve them separately and sum the obtained solutions, so the common solution will be found.
\begin{equation}\label{sec_mom_f_eq}
\begin{cases}
\displaystyle\frac {\partial{m_2(t,u)}}{\partial{t}} &= 2 m_2(t,u)\bigl(\beta - \mu\bigr)\enspace,\\[4pt]
m_2 (0,u) &= (En(0,u))^2 (1 - \delta_0 (u)) + \delta_0 (u) En^2(0,u)\enspace;
\end{cases}
\end{equation}

\begin{equation}\label{sec_mom_sec_eq}
\begin{cases}
\displaystyle\frac {\partial{m_2(t,u)}}{\partial{t}} &= 2 m_2(t,u)\bigl(\beta - \mu\bigr) + \frac {2k^2}{\mu - \beta}\enspace,\\[4pt]
m_2 (0,u) &= 0\enspace;
\end{cases}
\end{equation}

\begin{equation}\label{sec_mom_th_eq}
\begin{cases}
\displaystyle\frac {\partial{m_2(t,u)}}{\partial{t}} &= 2 m_2(t,u)\bigl(\beta - \mu\bigr) + 2\varkappa \mathcal L_{au} m_2 (t,u) - \frac{2\varkappa k a(u)}{\mu - \beta}\\&\quad + \delta_0 (u) \frac{k(2\mu + \sum_{n \geqslant 2} (n-1)(n-2) b_n)}{\mu - \beta}\enspace,\\[4pt]
m_2 (0,u) &= 0\enspace.
\end{cases}
\end{equation}

At the beginning we solve the first equation~\eqref{sec_mom_f_eq}. To reach this aim we are going to use Feinman-Kac formula (see~\cite{Ok05:e}). Then the solution $m_{2,1} (t,u)$ of~(\ref{sec_mom_f_eq}) is
\begin{align}\label{sec_mom_sol_1}
m_{2,1} (t,u) &= E[e^{- \int_0^{t} -2(\beta - \mu)ds} ((En(0,u))^2 (1 - \delta_0 (u)) + \delta_0 (u) En^2(0,u))]=\notag\\& = e^{2(\beta - \mu)t}E[(En(0,u))^2 (1 - \delta_0 (u)) + \delta_0 (u) En^2(0,u)]\enspace.
\end{align}

Note that this solution tends to zero as $t \to \infty$.

The second equation~(\ref{sec_mom_sec_eq}) is the ordinary differentional equation and the solution $m_{2,2} (t,u)$ can be found:
\begin{align}\label{sec_mom_sol_2}
m_{2,2} (t,u) = \frac{k^2}{(\mu - \beta)^2}(1-e^{2(\beta-\mu)t})\enspace.
\end{align}

As $t \to \infty$,
\begin{align*}
m_{2,2} (t,u) \to \frac{k^2}{(\mu-\beta)^2}\enspace.
\end{align*}

Finally, find the solution of the last equation~(\ref{sec_mom_th_eq}) which we denote by $m_{2,3} (t,u)$. Here we will use the discrete Fourier transform defined as
\begin{align}\label{dFtr}
\widehat f(\theta) = \sum_{u\in\mathbb Z^d} e^{i(\theta,u)} f(\theta),\qquad \theta \in [-\pi,\pi]^d\enspace.
\end{align}

In the last equation there is a term
\begin{align*}
\mathcal L_{au} m_{2,3} (t,u) &= \sum_{z \ne 0} a(z)\bigl(m_{2,3} (t,u+z) - m_{2,3} (t,u)\bigr)\\ &= \sum_{z \ne 0} a(z) m_{2,3} (t,u-z) - m_{2,3} (t,u)\enspace.
\end{align*}

The first term is the convolution of the functions $a(z)$ and $m_{2,3} (t,u)$. So applying the discrete Fourier transform~(\ref{dFtr}) to this term shows that
\begin{align*}
\widehat {\mathcal L_{au} m_{2,3}} (t,\theta) = \widehat {a}(\theta) \widehat {m_{2,3}} (t,\theta) - \widehat {m_{2,3}} (t,\theta)\enspace.
\end{align*}

Turn into discrete Fourier transform~(\ref{dFtr}) in the third differentional equation
\begin{equation*}
\begin{cases}
\displaystyle\frac {\partial \widehat m_{2,3} (t,\theta)}{\partial t} &= 2 \widehat m_{2,3} (t,\theta)[\beta - \mu] + 2 \varkappa (\widehat a(\theta) - 1)\widehat m_{2,3} (t,\theta)\\[4pt]&\quad - \frac {2 \varkappa k \widehat a(\theta)}{\mu - \beta} + \frac{k(2\mu + \sum_{n \geqslant 2} (n-1)(n-2) b_n)}{\mu - \beta}\enspace,\\[4pt]
\widehat{m_{2,3}}(\theta,0) &= 0\enspace.
\end{cases}
\end{equation*}

The solution of this equation has the form
\begin{align*}
\widehat m_{2,3}(t,\theta) = \frac{-\frac {2 \varkappa k \widehat a(\theta)}{\mu - \beta} + \frac{k(2\mu + \sum_{n \geqslant 2} (n-1)(n-2) b_n)}{\mu - \beta}}{2(\beta - \mu) + 2 \varkappa (\widehat a(\theta) - 1)} \Bigl(e^{(2(\beta - \mu) + 2 \varkappa (\widehat a(\theta) - 1))t} - 1\Bigr)\enspace.
\end{align*}

As $t \to \infty$:
\begin{align*}
\widehat m_{2,3}(t,\theta) \to \widehat m_{2,3}(\theta) = - \frac{-\frac {2 \varkappa k \widehat a(\theta)}{\mu - \beta} + \frac{k(2\mu + \sum_{n \geqslant 2} (n-1)(n-2) b_n)}{\mu - \beta}}{2(\beta - \mu) + 2 \varkappa (\widehat a(\theta) - 1)}\enspace.
\end{align*}

Then
\begin{align}\label{th_sol_F}
-\widehat m_{2,3}(\theta) = - \frac{C_1 \widehat a(\theta) + C_2}{C_3 - C_4 \widehat a(\theta)}\enspace,
\end{align}
where $C_1 = \frac{k\varkappa}{\mu-\beta}$, $C_2 = \frac{-k(\mu+\sum_{n\geqslant2} \frac{(n-1)(n-2)}{2}b_n)}{\mu-\beta}$, $C_3 = \mu-\beta+\varkappa$, $C_4 = \varkappa$.

From~(\ref{th_sol_F}) we receive
\begin{align*}
-\widehat m_{2,3}(\theta) = \frac{C_1}{C_4} + \left(\frac{C_1C_3 + C_2C_4}{C_4^2}\right) \frac{1}{\frac{C_3}{C_4}-\widehat a(\theta)}\enspace.
\end{align*}

\begin{remark}
The following properties are true
\begin{enumerate}
\item $\widehat \delta_0 (\theta) = 1$;
\item $\frac{1}{\frac{C_3}{C_4}-\widehat a(\theta)} = \frac{C_4}{C_3} \frac{1}{1-\frac{C_4}{C_3}\widehat a(\theta)} = \frac{C_4}{C_3} \left(1 + \sum_{n=1}^{\infty} \left(\frac{C_4}{C_3}\right)^n [\widehat a(\theta)]^n\right)$, if $|\frac{C_4}{C_3}\widehat a(\theta)| < 1$.
\end{enumerate}
\end{remark}

Then
\begin{align*}
-m_{2,3}(u) = -\frac{C_1}{C_4} \delta_0 (u) + \left(\frac{C_1}{C_4} + \frac{C_2}{C_3}\right) \left(\delta_0(u) + \sum_{n=1}^{\infty} \left(\frac{C_4}{C_3}\right)^n a^{*(n)} (u) \right)\enspace.
\end{align*}

Here $a^{*(n)} (u)$ - $n$-convolution of the function $a(z)$, thus is $a^{*(n)} (u) = (\underbrace{a*\ldots*a}_{n})(u)$.

After substituting the values of coefficients $C_1$, $C_2$, $C_3$, $C_4$ we get the following result:
\begin{align}\label{sec_mom_sol_3}
m_{2,3}(u) &= \frac{k}{\mu-\beta}\delta_0(u) + \frac{k}{(\mu-\beta)(\mu-\beta+\varkappa)} \left( \sum_{n=2}^{\infty} {n \choose 2} b_n  - \varkappa \right)\notag\\&\times \left( \delta_0(u) + \sum_{n=1}^{\infty} \left(\frac{\varkappa}{\mu-\beta+\varkappa}\right)^n a^{*(n)} (u) \right)\enspace.
\end{align}

Gather results obtained in~(\ref{sec_mom_sol_1}),~(\ref{sec_mom_sol_2}),~(\ref{sec_mom_sol_3})
\begin{align*}
m_2(t,u) &= m_{2,1} (t,u) + m_{2,2} (t,u) + m_{2,3} (t,u)\enspace.
\end{align*}

Then
\begin{align*}
m_2(t,u) \to \frac{k^2}{(\mu-\beta)^2} &+ \frac{k}{\mu-\beta}\delta_0(u) + \frac{k}{(\mu-\beta)(\mu-\beta+\varkappa)} \left( \sum_{n=2}^{\infty} {n \choose 2} b_n  - \varkappa \right) \\&\times \left( \delta_0(u) + \sum_{n=1}^{\infty} \left(\frac{\varkappa}{\mu-\beta+\varkappa}\right)^n a^{*(n)} (u) \right),\quad t \to \infty\enspace.
\end{align*}

\section{The Higher Moments}
\label{n_moment}

\quad Give the definition of the $n^{th}$ moment:
\begin{align*}
m_n(t, x_1,..., x_n) = E\Bigl[\prod_{i=1}^{n} n(t,x_i)\Bigr]\enspace.
\end{align*}

Let $\mathbf 1_{A}$ is the indicator of the set $A$, thus is
\begin{equation*}
\mathbf 1_{A} = \begin{cases}
1, &\text{if $A$ is true}\enspace;\\
0, &\text{otherwise}\enspace.
\end{cases}
\end{equation*}

Using the same methods of calculations (and relations~\eqref{first_ce}-\eqref{eighth_ce},~\eqref{indep} and~\eqref{fullprob}) it is possible to receive the differentional equations for higher moments.

\begin{align*}
m_n(t&+dt, x_1,...,x_n) = \Expect\Bigl[\prod_{i=1}^{n} n(t+dt,x_i)\Bigr] = \Expect\Bigl[\Expect\Bigl[\prod_{i=1}^{n} n(t+dt,x_i)\Bigr] |\mathcal F_{\le t}\Bigr] \\&= \Expect\Bigl[\Expect\Bigl[\prod_{i=1}^{n} (n(t,x_i) + \xi(dt,x_i))\Bigr]| \mathcal F_{\le t}\Bigr] = m_n(t,x_1,...,x_n) \\&+ \sum_{i=1}^{n}\Expect\Bigl[\prod_{j=1, j \ne i}^{n} [n(t,x_j)\Expect[[\xi(dt,x_i)]| \mathcal F_{\le t}]]\Bigr] \\&+ \sum_{i=1}^{\infty} \sum_{p=2}^{n}\Expect\Bigl[\prod_{x_j \ne x_i} [n(t,x_j) \Expect[[(\xi(dt,x_i))^p]| \mathcal F_{\le t}]]\Bigr] \mathbf 1_{A(p,x_i)} \\&+ \sum_{\substack{i,j=1,\\x_i \ne x_j}}^{n} \sum_{\substack{p,q>0:\\2 \le p+q \le n}} \Expect\Bigl[\prod_{x_k: x_k \ne x_i,x_j}[n(t,x_k)\Expect[[(\xi(dt,x_i))^p(\xi(dt,x_j))^q] \mathcal F_{\le t}]]\Bigr] \mathbf 1_{A(p,x_i)} \mathbf 1_{A(q,x_j)}\enspace.
\end{align*}

To continue calculations we use equalities~(\ref{first_ce})---(\ref{eighth_ce}):
\begin{align*}
m_n(t&+dt,x_1,...,x_n) = m_n(t,x_1,...,x_n) + \sum_{i=1}^{n} \Expect[\prod_{j=1, j \ne i}^{n} [n(t,x_j)[\sum_{r=2}^{\infty} (r-1)b_rn(t,x_i)dt + kdt \\&+ \sum_{z \ne 0} a(z)n(t,x_i+z)dt - (\mu + \varkappa)n(t,x_i)dt]]] + \sum_{i=1}^{n} \sum_{p=2}^{n} \Expect[\prod_{x_j \ne x_i} n(t,x_j) \\&\times [\sum_{r=2}^{\infty} (r-1)^p b_rn(t,x_i)dt + kdt + \varkappa \sum_{z \ne 0} a(z)n(t,x_i+z)dt + (-1)^p(\mu + \varkappa) \\&\times n(t,x_i)dt]] \mathbf 1_{A(p,x_i)} + \sum_{\substack{i,j=1,\\x_i \ne x_j}}^{n} \sum_{\substack{p,q>0:\\2 \le p+q \le n}} \Expect[\prod_{x_k:x_k \ne x_i,x_j} n(t,x_k)[(-1)^p \varkappa a(x_j - x_i)n(t,x_i) \\&\times dt + (-1)^q \varkappa a(x_i - x_j)n(t,x_j)dt]] \mathbf 1_{A(p,x_i)} \mathbf 1_{A(q,x_j)} = m_n(t,x_1,...,x_n) + n[\beta - \mu] \\&\times m_n(t,x_1,..,x_n)dt  + k\sum_{i=1}^{n} m_{n-1} (t,x_1,...,\hat x_i,...,x_n)dt + \varkappa \sum_{i=1}^{n} \mathcal L_{ax_i} m_n(t,x_1,...,x_n) \\&\times dt + \sum_{i=1}^{n} \sum_{p=2}^{n}[m_{n-p+1} (t,x_i,\tilde x_1,...,\tilde x_{n-p})[\sum_{r=2}^{\infty} (r-1)^p b_r + (-1)^p \mu]dt \\&+ km_{n-p} (t,\tilde x_1,...,\tilde x_{n-p})dt + \varkappa \mathcal L_{ax_i} m_{n-p+1} (t,x_i,\tilde x_1,...,\tilde x_{n-p})dt + (1 + (-1)^p)\varkappa \\&\times m_{n-p+1} (t,x_i,\tilde x_1,...,\tilde x_{n-p})dt] \mathbf 1_{A(p,x_i)} + \sum_{x_i \ne x_j} \sum_{\substack{p,q>0:\\2 \le p+q \le n}} [m_{n-(p+q) +1} (t,x_i,\tilde x_1,...,\tilde x_{n-(p+q)})\\&\times(-1)^p\varkappa a(x_j - x_i)dt + m_{n-(p+q)+1} (t,x_j,\tilde x_1,...,\tilde x_{n-(p+q)})(-1)^q\varkappa a(x_i - x_j)dt]\\&\times \mathbf 1_{A(p,x_i)} \mathbf 1_{A(q,x_j)} \enspace,
\end{align*}
where $A(p,x_i)$ is the following event:
\begin{align*}
A(p,x_i) = \{ \text{sequence } \{ x_1,...x_n \} \text{ contains exactly } p \text{ identical elements } x_i \}.
\end{align*}

Hence, from the representation above, we get the differentional equation
\begin{equation*}
\begin{cases}
\displaystyle\frac{\partial m_n(t,x_1,...,x_n)}{\partial t} &= n[\beta - \mu]m_n(t,x_1,..,x_n)  + k\sum_{i=1}^{n} m_{n-1} (t,x_1,...,\hat x_i,...,x_n) \\&+ \varkappa \sum_{i=1}^{n} \mathcal L_{ax_i} m_n(t,x_1,...,x_n) + \sum_{i=1}^{n} \sum_{p=2}^{n}[m_{n-p+1} (t,x_i,\tilde x_1,\\&...,\tilde x_{n-p})[\sum_{r=2}^{\infty} (r-1)^p b_r + (-1)^p \mu] + km_{n-p} (t,\tilde x_1,...,\tilde x_{n-p}) \\&+ \varkappa \mathcal L_{ax_i} m_{n-p+1} (t,x_i,\tilde x_1,...,\tilde x_{n-p}) + (1 + (-1)^p)\varkappa m_{n-p+1} (t,x_i,\tilde x_1,\\&...,\tilde x_{n-p})dt] \mathbf 1_{A(p,x_i)} + \sum_{x_i \ne x_j} \sum_{\substack{p,q>0:\\2 \le p+q \le n}} [m_{n-(p+q) +1} (t,x_i,\tilde x_1,\\&...,\tilde x_{n-(p+q)})(-1)^p\varkappa a(x_j - x_i) + m_{n-(p+q)+1} (t,x_j,\tilde x_1,...,\tilde x_{n-(p+q)}) \\&\times(-1)^q\varkappa a(x_i - x_j)] \mathbf 1_{A(p,x_i)} \mathbf 1_{A(q,x_j)}\enspace, \\
m_n(0,x_1,...,x_n) &= \sum_{\substack{p_1,...,p_n:\\p_1+...+p_n=n}}\Expect\prod_{i=1}^{n} n(0,x_i)^{p_i} \prod_{i=1}^{n} \mathbf 1_{A(p_i,x_i)}\enspace.
\end{cases}
\end{equation*}

Here record $m_{n-p+1} (t,x_i,\tilde x_1,...,\tilde x_{n-p})$ means that when we sum by $i$ there are no points $x_i$ in the set $\{ \tilde x_1,...,\tilde x_{n-p} \}$.

\section{Generating Function}
\label{gen_function}

\quad In this section we consider one of the most popular tools in BRWs studies. This is generating function.

\subsection{Definition and the Differentional Equation of the Generating Function.}
\label{def_eq_gen_func}

\quad Our target is to obtain the differentional equation of the generating function to simplify calculations in the moments' research.

\subsubsection{Definition.}
\label{def_gen_func}

\quad Define the generating function $F_{\infty} (z,t,x)$ of the particle field $n(t,x)$ as
\begin{equation*}
\begin{cases}
F_{\infty} (z,t,x) = Ee^{-zn(t,x)}\\
F_{\infty} (z,0,x) = Ee^{-zn(0,x)}\enspace
\end{cases},\quad z \geqslant 0\enspace.
\end{equation*}

Generating function is used in calculating moments (to do it, we should take the derivative of the variable $z$ and substitute $z=0$). Here we want to derivate the differentional equation for the generating function and compare with the results above.

\subsubsection{Differentional Equation of $F_{\infty} (z,t,x)$.}
\label{eq_gen_func}

\quad Derivating of the equation for the function $F_{\infty} (z,t,x)$ is based on the considering this function at the time $t+dt$. Then using~\eqref{nineth_ce}

\begin{align*}
F_{\infty} (z,t+dt,x) &= \Expect e^{-zn(t+dt,x)} = \Expect e^{-zn(t,x)}e^{-z\xi(dt,x)} =^{(2)} \Expect\Bigl(\Expect(e^{-zn(t,x)}e^{-z\xi(dt,x)}|\mathcal F_{\leqslant t})\Bigr) \\&= \Expect\Bigl(e^{-zn(t,x)}\Bigl(\Expect(e^{-z\xi(dt,x)}|\mathcal F_{\leqslant t})\Bigr)\Bigr)\enspace.
\end{align*}

Use~(\ref{nineth_ce}). Therefore,
\begin{align*}
F_{\infty} (z,t+dt,x) &= \Expect(e^{-zn(t,x)}(\sum_{n\geqslant2} e^{-z(n-1)}b_n n(t,x) dt + e^{-z}(kdt + \varkappa \sum_{y\ne0} a(y)n(t,x+y)dt) \\&+ e^z (\mu + \varkappa)n(t,x)dt + (1 - \sum_{n\geqslant2} b_n n(t,x)dt - kdt - \varkappa \sum_{y\ne0} a(y) n(t,x+y)dt \\&- \mu n(t,x)dt - \varkappa n(t,x)dt))) = F_{\infty} (z,t,x) + dt\Expect(e^{-zn(t,x)}(\sum_{n\geqslant2} b_n n(t,x) \\&+ ke^{-z} + \varkappa e^{-z} \sum_{y\ne0}a(y) n(t,x+y) + e^z (\mu + \varkappa)n(t,x) - \sum_{n\geqslant2} b_n n(t,x) \\&- k - \varkappa \sum_{y\ne0} a(y)n(t,x+y) - \mu n(t,x) - \varkappa n(t,x)))\enspace.
\end{align*}

Consequently,
\begin{align*}
\displaystyle\frac{F_{\infty} (z,t+dt,x) - F_{\infty} (z,t,x)}{dt} &= E(e^{-zn(t,x)}(\sum_{n\geqslant2} (e^{-z(n-1)} - 1) b_n n(t,x) + (e^{-z} - 1)k \\&+ (e^z - 1)(\mu + \varkappa)n(t,x) + (e^{-z} - 1)\sum_{y\ne0} a(y)n(t,x+y)))\enspace.
\end{align*}

\begin{remark}
Note that
\begin{enumerate}
\item $\displaystyle\frac{F_{\infty} (z,t+dt,x) - F_{\infty} (z,t,x)}{dt} \xrightarrow[dt\to0] \displaystyle\frac{\partial F_{\infty} (z,t,x)}{\partial t}$;
\item $\Expect\Bigl(e^{-zn(t,x)}n(t,x)\Bigr) = -\displaystyle\frac{\partial F_{\infty} (z,t,x)}{\partial z}$.
\end{enumerate}
\end{remark}

Thus, we receive the differentional equation of the generating function:
\begin{equation*}
\begin{cases}
\displaystyle\frac{\partial F_{\infty} (z,t,x)}{\partial t} &= \displaystyle\frac{\partial F_{\infty} (z,t,x)}{\partial z}\Bigl(\sum_{n\geqslant2} (1 - e^{-z(n-1)}) b_n + (1 - e^z)(\mu + \varkappa)\Bigr) + kF_{\infty} (z,t,x) \\&\times(e^{-z} - 1) + E\Bigl(e^{-zn(t,x)}(e^{-z} - 1)\sum_{y\ne0} a(y)n(t,x+y)\Bigr)\enspace,\\
F_{\infty} (z,0,x) &= Ee^{-zn(0,x)}\enspace.
\end{cases}
\end{equation*}

\subsection{Usage of the generating function}

\quad As mentioned above, generating functions help to obtain the equations for the moments of random variables.

Concretely,
\begin{align*}
m_1(t,x) &= \left. -\displaystyle\frac{\partial F_{\infty} (z,t,x)}{\partial t} \right|_{z=0}\enspace,\\
m_2(t,x,x) &= \left. \displaystyle\frac{\partial^2 F_{\infty} (t,x,y)}{\partial t^2} \right|_{z=0}\enspace.
\end{align*}

Now check results obtained for differentional equations of the first two moments in sections~\ref{first_moment} and~\ref{second_moment}.

\subsubsection{The First Moment}
\label{gen_func_f_mom}

\quad The equation for the first moment can be derived from
\begin{align*}
\left. \displaystyle\frac{\partial^2 F_{\infty} (z,t,x)}{\partial z\partial t} \right|_{z=0} = - \displaystyle\frac{\partial m_1(t,x)}{\partial t}\enspace.
\end{align*}

From this and~\ref{eq_gen_func} we get the left part
\begin{align*}
&\left. \displaystyle\frac{\partial F_{\infty} (z,t,x)}{\partial z} \right|_{z=0} \Bigl(\sum_{n\geqslant2} (n-1)b_n - \mu - \varkappa\Bigr) + kF_{\infty}(0,t,x) (-1) \\&+ \varkappa \Expect\sum_{y\ne0} a(y)n(t,x+y)(-1) = -m_1(t,x)(\beta-\mu) + \varkappa m_1(t,x) \\&- k - \varkappa \sum_{y\ne0} a(y)m_1(t,x+y)\enspace.
\end{align*}

The initial condition is
\begin{align*}
m_1(0,x) = En(0,x)\enspace.
\end{align*}

Hence, combining the above results, there is Cauchy problem for the first moment
\begin{equation*}
\begin{cases}
\displaystyle\frac{\partial m_1(t,x)}{\partial t} &= (\beta - \mu)m_1(t,x) + \varkappa \mathcal L_a m_1(t,x) + k\enspace,\\
m_1(0,x) &= En(0,x)\enspace.
\end{cases}
\end{equation*}

\subsubsection{The Second Moment}
\label{gen_func_s_mom}

\quad Identically do for the second moment:
\begin{align*}
\left. \displaystyle\frac{\partial^3 F_{\infty} (z,t,x)}{\partial z^2 \partial t} \right|_{z=0} = m_2(t,x,x)\enspace.
\end{align*}

So, in the left part we have
\begin{align*}
&2\left. \displaystyle\frac{\partial^2 F_{\infty} (z,t,x)}{\partial z^2} \right|_{z=0} \Bigl(\sum_{n\geqslant2}(n-1)b_n - \mu - \varkappa\Bigr) + \left. \displaystyle\frac{\partial F_{\infty} (z,t,x)}{\partial z} \right|_{z=0} (-\sum_{n\geqslant2} (n-1)^2 b_n \\&- \mu - \varkappa) + 2k\left. \displaystyle\frac{\partial F_{\infty} (z,t,x)}{\partial z} \right|_{z=0} (-1) + kF_{\infty} (0,t,x) + \varkappa \Expect(\sum_{y\ne0} a(y)n(t,x+y) \\&\times(2e^{-zn(t,x)}n(t,x) + e^{-zn(t,x)}))|_{z=0} = 2(\beta-\mu)m_2(t,x,x) - 2\varkappa m_2(t,x,x) \\&+ m_1(t,x)(\sum_{n\geqslant2} (n-1)^2 b_n + \mu + \varkappa) + 2km_1(t,x) + k + 2\varkappa\sum_{y\ne0} a(y)m_2(t,x,x+y) \\&+ \varkappa\sum_{y\ne0} a(y)m_1(t,x+y) = 2(\beta - \mu)m_2(t,x,x) + 2\varkappa\mathcal L_{ax} m_2(t,x,x) + k + m_1(t,x) \\&\times(\sum_{n\geqslant2} (n-1)^2 b_n + \mu + 2\varkappa + 2k) + \varkappa\mathcal L_a m_1(t,x)\enspace.
\end{align*}

The initial condition:
\begin{align*}
m_2(t,x,x) = En^2(0,x)\enspace.
\end{align*}

Thus,
\begin{equation*}
\begin{cases}
\displaystyle\frac{\partial m_2(t,x,x)}{\partial t} &= 2(\beta - \mu)m_2(t,x,x) + 2\varkappa\mathcal L_{ax} m_2(t,x,x) + k \\&+ m_1(t,x)(\sum_{n\geqslant2} (n-1)^2 b_n + \mu + 2\varkappa + 2k) + \varkappa\mathcal L_a m_1(t,x)\enspace,\\
m_2(0,x,x) &= En^2(0,x)\enspace.
\end{cases}
\end{equation*}

\subsection{Summary}
\label{gen_func_sum}

\quad As we see in~\ref{first_moment},~\ref{eq_sec_mom} and~\ref{gen_func_f_mom},~\ref{gen_func_s_mom}, obtained equations are identical. therefore, using the generating function makes the derivating of the equations much easier. However it is not possible to receive equation for high-order moments; because taking the derivative of the variable $z$ and substituting $z = 0$ gives the correlation functions only in case when points (which used as arguments) are identical. But in the right side of the equations there are various configurations of points. So these equations cannot be solved.

\section{Lyapunov Stability}
\label{lyapunov_stab}

\quad In sections~\ref{first_moment} and~\ref{second_moment} there were obtained asymptotics for the first two moments in case of constant rates $\beta$, $\mu$, $k$. In this part we are going to explore moments when the intencities are functions which depend on the point on the lattice, thus is $\beta = \beta(x)$, $\mu = \mu(x)$, $k = k(x)$. $x \in \mathbb Z^d$. For the next study there will be used the following

\begin{definition}~\cite[ch.~2,\S~10]{demidovich}
It is said that two systems of differentional equations
\begin{align*}
\frac{dx}{dt} = f(t,x)\quad and \quad \frac{dy}{dt} = g(t,y)
\end{align*}
are \textbf{asymptotically equivalent} systems if a one-to-one correspondence can be estableshed between their solutions $x(t)$, $y(t)$ respectively:
\begin{align*}
\lim_{t\to\infty} [x(t) - y(t)] = 0\enspace.
\end{align*}
\end{definition}

\subsection{Parabolic Problem}
\label{parab_prob}

\quad To study the Lyapunov stability we are going to use the following

\begin{lemma}\label{lemma_parab_prob}
Consider the parabolic problem
\begin{align*}
\begin{cases}
\displaystyle \frac{\partial u(t,x)}{\partial t} &= \mathcal L_a u(t,x) - v(x)u(t,x) + f(t,x)\enspace,\\
u(0,x) &= u_0(x),\quad x \in \mathbb Z^d,\quad t \in [0,T]\enspace.
\end{cases}
\end{align*}

Let $v(x)$, $u_0(x)$, $f(t,x)$ are bounded (the first two - on $\mathbb Z^d$, $f(t,x)$ - on $[0,T] \times \mathbb Z^d$), $x(t)$ is a symmetric random walk on $\mathbb Z^d$ with the generator $\mathcal L_a$. Then
\begin{align*}
u(t,x) = \Expect_x e^{- \int_0^t v(x(s))ds} u_0(x(t)) + \Expect_x \int_0^t f(t-s,x(s)) e^{- \int_0^s v(x(\tau))d\tau} ds\enspace.
\end{align*}
\end{lemma}

The proof of the Lemma~\ref{lemma_parab_prob} can be found in~\cite{ABMY}.

\subsection{Assumptions}
\label{assump}

\quad For the next research we make some assumptions. Let $k_0$, $\beta_0$, $\mu_0$, $u_0$, $\tilde u_0 \in \mathbb R$  and for fixed $\varepsilon > 0$ say
\begin{subequations}
\begin{equation} \label{a1}
\mu_0 - \beta_0 = v_0 > 0\enspace;
\end{equation}
\begin{equation}\label{a2}
k_0 > 0\enspace;
\end{equation}
\begin{equation}\label{a3}
u_0 > 0\enspace;
\end{equation}
\begin{equation}\label{a4}
\tilde u_0 > 0\enspace;
\end{equation}
\begin{equation}\label{a5}
v_0 - \varepsilon \leqslant v(x) \leqslant v_0 + \varepsilon\enspace;
\end{equation}
\begin{equation}\label{a6}
k_0 - \varepsilon \leqslant k(x) \leqslant k_0 + \varepsilon\enspace;
\end{equation}
\begin{equation}\label{a7}
u_0 - \varepsilon \leqslant u_0(x) \leqslant u_0 + \varepsilon\enspace;
\end{equation}
\begin{equation}\label{a8}
\tilde u_0 - \varepsilon \leqslant u_0(x,y) \leqslant \tilde u_0 + \varepsilon\enspace.
\end{equation}
\end{subequations}

\subsection{The First Moment}
\label{l_s_f_mom}

\quad Consider the equation for the first moment in case of non-constant rates:
\begin{align*}
\begin{cases}
\displaystyle \frac{\partial m_1(t,x)}{\partial t} &= \mathcal L_a m_1(t,x) - v(x)m_1(t,x) + k(x)\enspace,\\
m_1(0,x) &= u_0(x)\enspace.
\end{cases}
\end{align*}

\begin{theorem}
Under the assumptions in~\ref{assump} and for any $\medspace t \geqslant 0$
\begin{align*}
C_1^- \varepsilon + C_0^- e^{-(v_0+\varepsilon)t} \leqslant m_1(t,x) - \frac{k_0}{v_0} \leqslant C_1^+ \varepsilon + C_0^+ e^{-(v_0-\varepsilon)t}\enspace.
\end{align*}
Constants $C_0^{\pm}$, $C_1^{\pm}$ only depend on $k_0$, $v_0$, $u_0$ if $\varepsilon \leqslant \frac{\min(k_0,v_0)}{2}$.
\end{theorem}
\begin{remark}
$\frac{k_0}{v_0} = \frac{k_0}{\mu_0-\beta_0}$ is the limit vlue of the first moment in case of constant intencities when $t \to \infty$.
\end{remark}
\textbf{Proof:}
\begin{enumerate}
\item For the upper estimate we use assumptions~(\ref{a5}),~(\ref{a6}) and~(\ref{a7}):
\begin{align*}
m_1(t,x) &\leqslant (u_0 + \varepsilon) e^{-(v_0-\varepsilon)t} + \int_0^t (k_0+\varepsilon) e^{-(v_0-\varepsilon)s} ds \\&= (u_0 + \varepsilon) e^{-(v_0-\varepsilon)t} + \frac{k_0+\varepsilon}{v_0-\varepsilon}(1-e^{-(v_0-\varepsilon)t}) \\&= \frac{k_0+\varepsilon}{v_0-\varepsilon} + e^{-(v_0-\varepsilon)t} (u_0 - \frac{k_0+\varepsilon}{v_0-\varepsilon}) \\&= \frac{k_0}{v_0} + O(\varepsilon) + e^{-(v_0-\varepsilon)t} (u_0 - \frac{k_0}{v_0} + O(\varepsilon))\enspace.
\end{align*}
\item Use the same technique the lower estimate:
\begin{align*}
m_1(t,x) &\geqslant (u_0 - \varepsilon) e^{-(v_0+\varepsilon)t} + \int_0^t (k_0-\varepsilon) e^{-(v_0+\varepsilon)s} ds \\&= (u_0 - \varepsilon) e^{-(v_0+\varepsilon)t} + \frac{k_0-\varepsilon}{v_0+\varepsilon}(1-e^{-(v_0+\varepsilon)t}) \\&= \frac{k_0-\varepsilon}{v_0+\varepsilon} + e^{-(v_0+\varepsilon)t} (u_0 - \frac{k_0-\varepsilon}{v_0+\varepsilon}) \\&= \frac{k_0}{v_0} + O(\varepsilon) + e^{-(v_0+\varepsilon)t} (u_0 - \frac{k_0}{v_0} + O(\varepsilon))\enspace.
\end{align*}
\end{enumerate}
Now there is the assertion of the theorem.$\blacksquare$

From this theorem we get that the solutions for the equations in cases of constant and non-constant rates are asymptoticaly equivalent.

\subsection{The Second Moment}
\label{l_s_s_mom}

\quad In the equation of the second moment there are two operators $\mathcal L_{ax}$, $\mathcal L_{ay}$, and the equation has the form:
\begin{align*}
\begin{cases}
\displaystyle \frac{\partial m_2(t,x,y)}{\partial t} &= \mathcal L_{ax} m_2(t,x,y) + \mathcal L_{ay} m_2(t,x,y) - (v(x) + v(y))m_2(t,x,y) \\&+ k(x)m_1(t,y) + k(y)m_1(t,x) - \varkappa(a(x-y)m_1(t,x) \\&+ a(y-x)m_1(t,y)) + \delta_x (y) (m_1(t,x) (\mu(x) + \sum_{n\geqslant2} (n-1)^2 b_n(x)) \\&+ k(x) + \varkappa\mathcal L_a m_1(t,x))\enspace,\\
m_2(0,x,y) &= u_0(x,y)\enspace.
\end{cases}
\end{align*}

Let
\begin{align}
F(t,x) = m_1(t,x) ( \mu(x) + \sum_{n=2}^{\infty} (n-1)^2 b_n(x)) + k(x) + \varkappa\mathcal L_a m_1(t,x)\enspace; \label{F_func_s_m}
\end{align}
\begin{align}
f(t,x,y) &= k(x)m_1(t,y) + k(y)m_1(t,x) + \delta_x(y) F(t,x) \notag\\&- \varkappa a(x-y)(m_1(t,y) + m_1(t,x))\enspace; \label{f_func_s_m}
\end{align}
\begin{align}
V(x,y) = v(x) + v(y)\enspace. \label{V_func_s_m}
\end{align}

Then the equation can be rewright:
\begin{align*}
\begin{cases}
\displaystyle \frac{\partial m_2(t,x,y)}{\partial t} &= (\mathcal L_{ax} + \mathcal L_{ay}) m_2(t,x,y) - V(x,y)m_2(t,x,y) + f(t,x,y)\enspace,\\
m_2(0,x,y) &= u_0(x,y)\enspace.
\end{cases}
\end{align*}

For this equation we also want to apply lemma~\ref{lemma_parab_prob}, so we reckon couple of two independent random walks on $\mathbb Z^d$: $x(t)$ and $y(t)$ with generators $\mathcal L_{ax}$ and $\mathcal L_{ay}$ respectively. Thus, now we can use lemma~\ref{lemma_parab_prob} for the pair $(x(t),y(t))$.

From this proposition we receive
\begin{align*}
m_2(t,x,y) &= \Expect_{(x,y)} e^{-\int_0^t V(x(s),y(s))ds} u_0(x(t),y(t)) \\&+ \Expect_{(x,y)} \int_0^t f(t-s,x(s),y(s))e^{-\int_0^s V(x(\tau),y(\tau))d\tau}ds\enspace.
\end{align*}

\begin{remark}
From~(\ref{a5}) and~(\ref{V_func_s_m})
\begin{align*}
2(v_0 - \varepsilon) \leqslant V(x,y) \leqslant 2(v_0 + \varepsilon)\enspace.
\end{align*}
\end{remark}

Let
\begin{align} \label{G_func_s_m}
G(t,x,y) := \Expect_{(x,y)} e^{-\int_0^t V(x(s),y(s))ds} u_0(x(t),y(t))\enspace;
\end{align}
\begin{align} \label{H_func_s_m}
H(t,x,y) := \Expect_{(x,y)} \int_0^t f(t-s,x(s),y(s))e^{-\int_0^s V(x(\tau),y(\tau))d\tau}ds\enspace.
\end{align}

So
\begin{align} \label{s_m_sum_G_H}
m_2(t,x,y) = G(t,x,y) + H(t,x,y)
\end{align}

\begin{enumerate}

\item Estimate function $G(t,x,y)$ (defined in~\eqref{G_func_s_m}) using~(\ref{a8}). Then
\begin{enumerate}
\item $G(t,x,y) \leqslant \Expect_{(x,y)} e^{-\int_0^t 2(v_0 - \varepsilon)ds}(u_0 + \varepsilon) = e^{-2(v_0 - \varepsilon)t}(\tilde u_0 + \varepsilon)$\enspace;
\item $G(t,x,y) \geqslant \Expect_{(x,y)} e^{-\int_0^t 2(v_0 + \varepsilon)ds}(u_0 - \varepsilon) = e^{-2(v_0 + \varepsilon)t}(\tilde u_0 - \varepsilon)$\enspace.
\end{enumerate}

\item Estimate $H(t,x,y)$ (from~\eqref{H_func_s_m}).
\begin{enumerate}
\item Firstly, find the boundaries for function $f(t,x,y)$ (definition see in~(\ref{f_func_s_m})).

\begin{remark}
\begin{equation*}
a(x-y) =
\begin{cases}
-1, &\text{if}~x=y\\
\in[0,1], &\text{if}~x\ne y\\
\end{cases}
\end{equation*}
$\Rightarrow$
\begin{equation*}
-a(x-y) =
\begin{cases}
1, &\text{if}~x=y\\
\in[-1,0], &\text{if}~x\ne y\\
\end{cases}
\end{equation*}
\end{remark}
\begin{enumerate}
\item Using assumptions~(\ref{a5}) and~(\ref{a6}) gives the upper boundary
\begin{align*}
f(t,x,y) &\leqslant (k_0 + \varepsilon + \varkappa)(m_1(t,x) + m_1(t,y)) + \delta_x(y) F(t,x) \\&= 2(k_0 + \varepsilon + \varkappa)\Bigl(\frac{k_0}{v_0} + C_1^+ \varepsilon + C_0^+ e^{-(v_0-\varepsilon)t}\Bigr) + \delta_x(y) F(t,x)\enspace;
\end{align*}
\item and the lower boundary
\begin{align*}
f(t,x,y) &\geqslant (k_0 - \varepsilon - \varkappa)(m_1(t,x) + m_1(t,y)) + \delta_x(y) F(t,x) =\\&= 2(k_0 - \varepsilon)\Bigl(\frac{k_0}{v_0} + C_1^- \varepsilon + C_0^- e^{-(v_0+\varepsilon)t}\Bigr) \\&- 2\varkappa\Bigl(\frac{k_0}{v_0} + C_1^+ \varepsilon + C_0^+ e^{-(v_0-\varepsilon)t}\Bigr) + \delta_x(y) F(t,x)\enspace.
\end{align*}
\end{enumerate}

\item Estimate $H(t,x,y)$
\begin{enumerate}
\item For upper estimate use i form 2.(a).
\begin{align*}
H(t,x,y) &\leqslant \int_0^t \delta_x(y)F(t-s,x)e^{-\int_0^t V(x(\tau),y(\tau))d\tau}ds \\&+ 2(k_0 + \varepsilon + \varkappa)\int_0^t\Bigl(\frac{k_0}{v_0} + C_1^+ \varepsilon + C_0^+ e^{-(v_0-\varepsilon)(t-s)}\Bigr)e^{-2(v_0-\varepsilon)s}ds\enspace.
\end{align*}
Let
\begin{align*}
L(t,x,y) = \int_0^t \delta_x(y)F(t-s,x)e^{-\int_0^s V(x(\tau),y(\tau))d\tau}ds\enspace.
\end{align*}

Then
\begin{align*}
H(t,x,y) &\leqslant L(t,x,y) + \frac{(k_0 + \varepsilon + \varkappa)(\frac{k_0}{v_0} + C_1^+\varepsilon)}{v_0-\varepsilon}(1 - e^{-2(v_0-\varepsilon)t}) \\&+ C_0^+ e^{-(v_0-\varepsilon)t}2(k_0 + \varepsilon + \varkappa) \frac{1}{v_0-\varepsilon}(1 - e^{-(v_0-\varepsilon)t}) \\&= L(t,x,y) + \frac{k_0(k_0 + \varkappa)}{v_0^2} + O(\varepsilon) + \frac{2C_0^+(k_0 + \varepsilon + \varkappa)}{v_0}e^{-(v_0-\varepsilon)t} \\&- e^{-2(v_0-\varepsilon)t}\Bigl(\frac{k_0(k_0 + \varkappa)}{v_0^2} + O(\varepsilon) + \frac{2C_0^+(k_0+\varkappa)}{v_0}\Bigr)\enspace.
\end{align*}
\item For the lower estimate use ii from 2.(a).
\begin{align*}
H(t,x,y) &\geqslant L(t,x,y) + 2(k_0 - \varepsilon)\int_0^t\Bigl(\frac{k_0}{v_0} + C_1^- \varepsilon + C_0^- e^{-(v_0+\varepsilon)(t-s)}\Bigr)e^{-2(v_0+\varepsilon)s}ds \\&- 2\varkappa\int_0^t\Bigl(\frac{k_0}{v_0} + C_1^+ \varepsilon + C_0^+ e^{-(v_0-\varepsilon)(t-s)}\Bigr)e^{-2(v_0-\varepsilon)s}ds \\&= L(t,x,y) + \frac{k_0-\varepsilon}{v_0+\varepsilon}\Bigl(\frac{k_0}{v_0}+ C_1^- \varepsilon\Bigr)(1 - e^{-2(v_0+\varepsilon)t}) \\&+ \frac{2C_0^-(k_0-\varepsilon)}{v_0+\varepsilon}e^{-(v_0+\varepsilon)t}(1 - e^{-(v_0+\varepsilon)t}) \\&- \varkappa\Bigl(\frac{k_0}{v_0} + C_1^+ \varepsilon\Bigr)\frac{1}{v_0-\varepsilon}(1 - e^{-2(v_0-\varepsilon)t}) \\&- \frac{2\varkappa C_0^+}{v_0-\varepsilon}e^{-(v_0-\varepsilon)t}(1 - e^{-(v_0-\varepsilon)t}) \\&= L(t,x,y) + \frac{k_0^2}{v_0^2} + O(\varepsilon) - e^{-2(v_0+\varepsilon)t}\Bigl(\frac{k_0^2}{v_0^2} + O(\varepsilon)\Bigr) \\&+ \frac{2C_0^- k_0}{v_0}e^{-(v_0+\varepsilon)t} - \Bigl(\frac{2C_0^- k_0}{v_0} + O(\varepsilon)\Bigr)e^{-2(v_0+\varepsilon)t} \\&- \frac{\varkappa k_0}{v_0^2} + \Bigl(\frac{\varkappa k_0}{v_0^2} + O(\varepsilon)\Bigr)e^{-2(v_0-\varepsilon)t} - \frac{2\varkappa C_0^+}{v_0}e^{-(v_0-\varepsilon)t} \\&+ \Bigl(\frac{2\varkappa C_0^+}{v_0} + O(\varepsilon)\Bigr)e^{-2(v_0-\varepsilon)t}\enspace.
\end{align*}
\end{enumerate}

\end{enumerate}
\item From~(\ref{s_m_sum_G_H}) and 1. и 2. there is an estimate fo the second moment in case of non-constant intencities:
\begin{enumerate}
\item Upper limit
\begin{align*}
m_2(t,x,y) &\leqslant e^{-2(v_0 - \varepsilon)t}(u_0 + \varepsilon) + L(t,x,y) + \frac{k_0(k_0 + \varkappa)}{v_0^2} \\&+ O(\varepsilon) + \frac{2C_0^+(k_0 + \varepsilon + \varkappa)}{v_0}e^{-(v_0-\varepsilon)t} \\&- e^{-2(v_0-\varepsilon)t}\Bigl(\frac{k_0(k_0 + \varkappa)}{v_0^2} + O(\varepsilon) + \frac{2C_0^+(k_0+\varkappa)}{v_0}\Bigr) \\&= L(t,x,y) + \frac{k_0^2}{v_0^2} + \frac{k_0 \varkappa}{v_0^2} + C_2^+ e^{-(v_0-\varepsilon)t} \\&+ \Bigl(u_0 - \frac{k_0+\varkappa}{v_0}\Bigl(2C_0^+ + \frac{k_0}{v_0}\Bigr) + O(\varepsilon)\Bigl) e^{-2(v_0-\varepsilon)t} + O(\varepsilon)\enspace.
\end{align*}
\item Lower limit
\begin{align*}
m_2(t,x,y) &\geqslant e^{-2(v_0 + \varepsilon)t}(u_0 - \varepsilon) + L(t,x,y) + \frac{k_0^2}{v_0^2} \\&+ O(\varepsilon) - e^{-2(v_0+\varepsilon)t}\Bigl(\frac{k_0^2}{v_0^2} + O(\varepsilon)\Bigr) + \frac{2C_0^- k_0}{v_0}e^{-(v_0+\varepsilon)t} \\&- \Bigl(\frac{2C_0^- k_0}{v_0} + O(\varepsilon)\Bigr)e^{-2(v_0+\varepsilon)t} - \frac{\varkappa k_0}{v_0^2} + \frac{\varkappa k_0}{v_0^2}e^{-2(v_0-\varepsilon)t} \\&- \frac{2\varkappa C_0^+}{v_0}e^{-(v_0-\varepsilon)t} + \frac{2\varkappa C_0^+}{v_0}e^{-2(v_0-\varepsilon)t}\enspace.
\end{align*}
\end{enumerate}

\end{enumerate}

Finally, we obtain the following boundaries fo the second moment:
\begin{align*}
A \leqslant m_2(t,x,y) - L(t,x,y) - \frac{k_0^2}{v_0^2} \leqslant B\enspace,
\end{align*}
where
\begin{align*}
A = C_2^- e^{-(v_0+\varepsilon)t} + C_3^- e^{-2(v_0+\varepsilon)t} + C_4^- \varepsilon - \frac{\varkappa k_0}{v_0^2}(1 - e^{-2(v_0-\varepsilon)t})\enspace;
\end{align*}
\begin{align*}
B = C_2^+ e^{-(v_0-\varepsilon)t} + C_3^+ e^{-2(v_0-\varepsilon)t} + C_4^+ \varepsilon + \frac{\varkappa k_0}{v_0^2}(1 - e^{-2(v_0-\varepsilon)t})\enspace.
\end{align*}

\textbf{Acknowledgements.} Yu.Makarova and E.Yarovaya were supported by the Russian Foundation for the Basic Research (RFBR), project No. 17-01-00468. S.Molchanov was supported by the Russian Science Foundation (RSF), project No. 17-11-01098.

\renewcommand{\refname}{References}

\end{document}